\documentclass[12pt,a4paper,reqno]{amsart}
\usepackage{preemble, relsize, enumitem, mathrsfs}

\usepackage{times,amsfonts,amsmath,amstext,amsbsy,amssymb,
  amsopn,amsthm,upref,eucal, amscd}
\usepackage[colorlinks=true]{hyperref}
\usepackage[T1]{fontenc}
\usepackage[backgroundcolor=white]{todonotes}
\usepackage{multirow}
\usepackage{color}
\usepackage{graphicx}
\usepackage{accents}

\date{\today}

\author{Stephen Cantrell}
\address{Mathematical Institute, University of St Andrews, North Haugh, St Andrews, KY16 9SS, Scotland}
\email{sjc33@st-andrews.ac.uk}

\date{\today. 
}

\title[A local limit theorem for singularity visits on translation surfaces]{A local limit theorem for singularity visits on translation surfaces}

\begin{document}

\begin{abstract}
Let $(X,\omega)$ be a translation surface with at least one singularity. We investigate the distribution of the number of singularity encounters along prime closed geodesics on $X$.  Writing $\ell(\gamma)$ for the geometric length of a closed geodesic and $m(\gamma)$ for the number of singularities it encounters, we establish a uniform counting local limit theorem for the values of $m(\gamma)$ sampled over $\ell(\gamma)<T$. The key technical step is to introduce and study a trace-class operator whose spectral data encodes the Fourier transform of our counting distribution. 
\end{abstract}

\maketitle

\section{Introduction}
 A translation surface is a compact Riemann surface $X$ equipped with a  (non-trivial)  holomorphic one-form $\omega$.
In other words, $X$ is a compact surface equipped with a metric that is flat except at a finite set $\mathcal{X} = \{x_1, \ldots, x_n\}$ of points, called singularities, that have cone angles $2\pi (k(x_i)+1)$, where $k(x_i) \in \mathbb N$, for $i=1, \ldots, n$. A saddle connection is a straight line segment in $X$ that starts and ends at singularities and does not visit any other singularities. There has been particular interest in obtaining asymptotics for the number of saddle connections of length less than $T$. Such results were initiated by Masur's work \cite{Mas90} (which established coarse quadratic growth) and includes the asymptotic results of Eskin-Masur \cite{EM01}, Eskin-Masur-Zorich  \cite{EMZ03} and many other beautiful works \cite{athreya2012distribution, NevoRuhrWeiss2020, athreya2023counting, robertson2023uniform}.  In this work, instead of studying individual saddle connections, we consider the asymptotic properties of closed saddle geodesics: locally length-minimising closed paths, up to cyclic reordering, that pass through singularities.

For such a geodesic $\gamma$ there are two natural and quite different notions of length.
The first is its geometric length $\ell(\gamma)$, measured in the flat metric.  The
second is the number $m(\gamma)$  of singularity encounters in one
period, counted with multiplicity.  Comparing combinatorial and geometric lengths is a central theme in various geometric contexts. In the classical setting of negatively curved surfaces (and more generally manifolds), counting
results comparing such lengths have been obtained by many authors, \cite{milnor, Lalley1989, schwartz.sharp, babillot.ledrappier,  Sharp2004, CantrellPollicott2022}.  The singularities on $(X,\omega)$ can be thought of as concentrated points of negative curvature and it is therefore very natural to try to obtain statistical results, analogous to those for negatively curved surfaces, for translation surfaces.

Cantrell-Pollicott \cite{CP25} recently established counting limit laws for
saddle-connection paths on translation surfaces (extending the work of Colognese-Pollicott \cite{ColP}).
To recall their result, let $\SCon$ denote the saddle connections of $(X,\omega)$. For $s \in \SCon$ let $i(s), t(s) \in \mathcal X$ denote the initial and terminal singularity of $s$.  A saddle connection path is a string $(s_1,\ldots, s_n)$ of  saddle connections forms a locally length minimising path in $(X,\omega)$ if $t(s_j) = i(s_{j+1})$ for each $j = 1, \ldots, n-1$ and the angle between $s_j$, $s_{j+1}$ is at least $\pi$ (see Section 4 of \cite{CP25}). The geometric length of $p = (s_1,\ldots, s_n)$ is its distance traversed in $(X,\omega)$: $\ell(p) = \ell(s_1) + \cdots + \ell(s_n)$. Whereas $m(p) = n$ is the number of saddle connections making up $p$. For $x \in \mathcal X$, let $\mathcal P_x$ denote the collection of saddle connection paths starting at $x$. Let $\nu_T$ denote the uniform counting measure on $\{ p \in \mathcal P_x : \ell(p) < T\}$ for each $T>0$.

\begin{theorem}[\cite{CP25} Theorem 1.2]\label{thm.sing}
There exist $\Lambda>0$ and $\sigma^2>0$ such that, for every singularity
$x\in\mathcal X$ and every $a<b$,
\[
\lim_{T\to\infty}\nu_{T}\left\{p \in\mathcal P_x :\frac{m(\gamma)-\Lambda T}{\sqrt T}\in[a,b]\right\} =\frac{1}{\sqrt{2\pi\sigma^2}}\int_a^b \exp\!\left(-\frac{t^2}{2\sigma^2}\right)\,dt.
\]
\end{theorem}
The original statement in \cite{CP25} has $m(p)+1$ instead of $m(p)$ (i.e. it counts the number of singularities visited) but subtracting one leaves the central limit theorem unchanged. 

\subsection*{Main result}
The purpose of this note is to obtain a corresponding local limit theorem but for counting over prime closed geodesics. We decided to pursue counting over closed geodesics as it feels more in keeping with the aforementioned results for negatively curved surfaces.
Let $\Prim$ be the set of oriented prime closed geodesics. These are closed, saddle connection paths, considered up to cyclic permutation,  that are not iterates of a shorter closed geodesic (see Definition \ref{def:prime}). We set
$\Prim(T)=\{\gamma\in\Prim:\ell(\gamma)<T\}$ and recall that, by Theorem C in \cite{ColP} we have the asymptotic
\[
\#\Prim(T) \sim \frac{e^{hT}}{hT} \ \ \ \text{ as }  \ \ T\to\infty
\]
where $h >0$ is a constant.  For $\gamma\in\Prim$, the integer $m(\gamma)$
counts the saddle connections in one prime period (equivalently the number of singularities visited along the closed geodesic).   For $T>0$  let $\mu_T$ denote the uniform measure on
$\Prim(T)$, so that
\[
\mu_T(E)=\frac{\#(E\cap\Prim(T))}{\#\Prim(T)}
\]
for $E \subset \Prim$.
Our main theorem is the following uniform local limit theorem.

\begin{theorem}\label{thm:main}
For every translation surface $(X,\omega)$ with at least one
singularity, there exist $\Lambda>0$ and $\sigma^2>0$ such that
\[
\lim_{T\to\infty}\sup_{n\in\N}\left|\sqrt T\, \mu_T\left\{\gamma\in\Prim :m(\gamma)=n\right\} -\frac{1}{\sqrt{2\pi\sigma^2}} \exp \, \left(-\frac{(n-\Lambda T)^2}{2\sigma^2T}\right) \right|=0.
\]
The constants $\Lambda$ and $\sigma^2$ are the same as in
Theorem \ref{thm.sing}.
\end{theorem}

This result shows that the probability of seeing any fixed integer number of
singularity encounters has a Gaussian asymptotic at the scale $T^{-1/2}$. We immediately deduce the following corollaries.

\begin{corollary}
Fix a  translation surface $(X,\omega)$ with at least one
singularity. Then, for every integer $q \ge 2$ and $ 0 \le l \le q-1$ we have that
\[
\frac{\#\{\gamma \in \Prim(T) : m(\gamma) = l \pmod q\}}{\#\Prim(T)} \longrightarrow \frac{1}{q}
\]
as $T\to\infty$. That is, $m(\gamma)$ is uniformly distributed modulo $q$ (for any $q$), when sampling $\gamma \in \Prim$ with respect to the sequence of uniform counting measures $\mu_T$.
\end{corollary}

This is easy to deduce from Theorem \ref{thm:main} (see for example Theorem 1.3 in \cite{DolgopyatHafouta2022}).
We also have the following result, which follows from  \cite[Theorem C]{ColP} and setting $T = n/\Lambda$ in Theorem \ref{thm:main}.

\begin{corollary}
Let  $(X,\omega)$ be a translation surface with at least one
singularity. Write  $\Lambda>0, \sigma^2 > 0$ for the constants in Theorem \ref{thm:main}. Then
\[
\#\left\{ \gamma \in \Prim : m(\gamma) = n \ \text{and } \ \ell(\gamma) < \frac{n}{\Lambda}\right\} \sim \, \frac{\Lambda^{3/2}}{h\sqrt{2\pi\sigma^2}} \, \frac{e^{hn/\Lambda}}{n^{3/2}}
\]
as $n\to \infty$ (through the integers).
Here $h$ is the exponential growth rate of $\#\Prim(T)$ as $T\to \infty$.
\end{corollary}
This provides a precise asymptotic for the number of prime closed geodesics $\gamma$ with $m(\gamma) = n$ (as long as we scale the geometric length appropriately).

\subsection*{Outline of the proof}
Following \cite{ColP, CP25} we encode geodesics by cyclic paths in a  countable directed graph
whose vertices are oriented saddle connections.  The results in \cite{CP25}
provides key structural results, such as connectivity properties, for this graph.  We first show, using geometric arguments, that
the period of the graph is always one (Proposition \ref{prop:aperiodic}) independent of $(X,\omega)$. This removes
possible congruence restrictions on $m(\gamma)$.  

The key step is then to further develop the functional analytic methods used in \cite{CP25}. We define a complex parameter family of linear operators $A_r$ on $\ell^2(\SCon)$ and prove that they are
\begin{enumerate}
\item trace class; and,
 \item their traces can be used to enumerate closed geodesics.
 \end{enumerate}
We then apply the appropriate positive operator theory and analytic perturbation theory to these operators to deduce that they have a spectral gap and an analytically varying, lead eigenvalue,
We then obtain non-lattice spectral estimates for this operator using work from \cite{CP25}.  Together with the usual Fourier inversion and smoothing tricks (as in the classical local limit theorems of Stone \cite{Stone65}, Lalley \cite{Lalley1987}, etc)  this gives precise counts for closed geodesics.  Along the way, we identify the mean and variance constants in Theorem \ref{thm:main} with those appearing in \cite{CP25} (i.e. in Theorem \ref{thm.sing}). This identification avoids the need to establish the positivity of the constants directly from the definition of the trace class operator $A_r$ which appears difficult to do.

Section \ref{sec:prelim} introduces the saddle connection coding and proves preliminary structural results about it. Section \ref{sec:spectral} develops the trace-class operator
and its spectral properties, while Section \ref{sec:fixed} establishes the
local limit theorem and passes to prime geodesics.  The appendix identifies
the constants with those  in \cite{CP25}.

\subsection*{AI disclosure}
This project was the author's first attempt at using AI assistance in research. The author had considered this local limit problem before and had a full outline for the proof but was struggling to find an appropriate Perron-Frobenius result to apply to the operator in Section \ref{sec:spectral}. The author described the outline of their proof to ChatGPT and asked it to `find an appropriate Perron-Frobenius theorem' to apply to the operator defined in Section \ref{sec:spectral}. After a short discussion, the author found an appropriate result (\cite[Theorem 1.5]{CWW20}) and using it, was able to complete the proof. AI was not used in the write-up of this work, nor to assist with the mathematical arguments (other than the search for a useful Perron-Frobenius theorem). Lastly, AI was used  for a few literature searches, was used to proof-read, and made a few suggestions on (non-mathematical) wording throughout this work.

\section{Preliminaries: geometry and coding}\label{sec:prelim}
\subsection{Translation surfaces}
We briefly recap some introductory material concerning translation surfaces, see \cite{Wright15} for more details.
Let $\mathcal X$ be the nonempty set of zeros of $\omega$ and $\SCon$ the countable set of oriented saddle connections. For $s\in\SCon$, let $i(s),t(s)\in\mathcal X$ denote its endpoints (initial and terminal respectively) and $\ell(s)>0$ its geometric length. We write
\[
\ell_*:=\inf_{s\in\SCon}\ell(s)>0
\]
for the length of shortest saddle connection.
A concatenation of saddle connections is locally length-minimising (i.e. a local geodesic) precisely when, at each singularity visited, the angle between the incoming and outgoing saddle connections  are at least $\pi$; see \cite[Definition 4.2]{CP25}.

Following Colognese-Pollicott \cite{ColP} we encode possible closed geodesics using an infinite graph.

\begin{definition} \label{def:paths}
Write $s_1 \rightsquigarrow s_2$ if $t(s_1)=i(s_2)$ and $s_1s_2$ is a local geodesic. We define a directed graph $\mathscr G_X$ with vertex set $\SCon$ and an edge joining $s_1$ to $s_2$ precisely when $s_1 \rightsquigarrow s_2$. A closed saddle path $w=(s_1,\ldots,s_n)$ is a path with $s_j\rightsquigarrow s_{j+1}$ for $1\leq j\leq n-1$ and  $s_n \rightsquigarrow s_1$. Its geometric length is $\ell(w)=\sum_{j=1}^n\ell(s_j)$ and its saddle connection count is $m(w)=n$. Two closed paths represent the same oriented closed geodesic if they differ by a cyclic rotation. The functions $\ell$ and $m$ extend naturally to closed geodesics.
\end{definition}

\begin{remark}\label{rem:paths}
Here we are using a different terminology to that used in \cite{CP25}. Indeed, in \cite{CP25} a closed saddle path denotes  a saddle connection path $p= (s_1,\ldots, s_n)$ that starts and ends at the same singularity (but may not correspond to a closed geodesic path, i.e. we have $t(s_n) = i(s_1)$ but we may not have $s_n \rightsquigarrow s_1$). 
\end{remark}

We define the graph period of $\mathscr G_X$ to be
\begin{equation*}
d(\mathscr G_X) =\gcd\{n\geq1:\mathscr G_X\text{ has an admissible closed path with }n\text{ edges}\}
\end{equation*}
where $\gcd$ denotes the greatest common divisor.

\begin{definition}\label{def:prime}
A closed geodesic $\gamma$  corresponding to a closed saddle connection path $(s_1,\ldots,s_n)$ is prime if it is not the iterate of any closed geodesic consisting of strictly  fewer than  $n$ saddle connections.
\end{definition}

We require the following three conditions established in Section 4 of \cite{CP25}.

\begin{proposition} \label{prop:geo}
For every $(X,\omega)$ as above:
\begin{enumerate}[label=(\roman*)]
\item For every $t>0$, $\sum_{s\in\SCon}e^{-t\ell(s)}<\infty$;
\item There is $C>0$ such that for any $s_1,s_2\in\SCon$ there is a saddle connection path beginning with $s_1$, ending with $s_2$, of length at most $\ell(s_1)+\ell(s_2)+C$; and, 
\item The $\mathbb Q$-vector space spanned by $\{\ell(s):s\in\SCon\}$ has infinite dimension. Likewise, the $\mathbb Q$-span of the lengths of closed geodesics  is infinite-dimensional.
\end{enumerate}
\end{proposition}
\begin{proof}
See \cite{CP25} for (i) and (ii). Point (iii) is essentially \cite[Proposition 4.8]{CP25} but some care is needed. Indeed, \cite[Proposition 4.8]{CP25} states that the set of lengths of saddle paths that start and end at the same singularity generates an infinite rank $\mathbb Q$ vector space. Such paths need not be closed as in Definition \ref{def:paths} (see Remark \ref{rem:paths}). However, by (ii) there are a finite collection of saddle connections $F \subset \SCon$ such that any saddle connection path can be extended into a closed path using saddle connections in $F$. From this observation and \cite[Proposition 4.8]{CP25}  point (iii) follows easily.
\end{proof}

The following result shows that the period of $\mathscr G_X$ is always $1$, that is, there are no arithmetic constraints on how many 
saddle connections make up a closed geodesic.

\begin{proposition}\label{prop:aperiodic}
For every translation surface with at least one singularity, $d(\mathscr G_X)=1$, i.e. $\mathscr G_X$ has period one. 
\end{proposition}

\begin{proof}
Let $d=d(\mathscr G_X)$ and fix a base vertex $s_0\in\SCon$. For each $s\in\SCon$, define $j(s)\in\mathbb Z/d\mathbb Z$ to be the number of edges in a directed path from $s_0$ to $s$, taken modulo $d$.
This is well defined because if two such paths have lengths $n_1,n_2$, that we can, using Proposition \ref{prop:geo} (ii),
append both paths by the same directed path from $s$ to $s_0$. The resulting
closed paths have lengths $n_1+k$ and $n_2+k$ which are both divisible by $d$.
It follows that 
\begin{equation}\label{eq:period-class-transition}
 j(s')=j(s)+1\pmod d \ \ \ \ \text{whenever }s\rightsquigarrow s'.
\end{equation}

We claim that $j(s)$ depends only on the initial singularity $i(s)$.
Fix a singularity $x \in \mathcal X$ and choose a triangulation of $X$ by saddle
connections,  \cite[Lemma 1.11]{Wright15}. Write $L(x)=2\pi(k(x)+1)\geq 4\pi$ for the total cone angle
at $x \in \mathcal X$. The outgoing directions of the triangulation edges at $x \in \mathcal X$ divide
the possible angles $\mathbb R/L(x)\mathbb Z$ into finitely many arcs,
each of length strictly less than $\pi$. In particular, every open angle arc
of length greater than $\pi$ contains an outgoing triangulation-edge angle.

Suppose first that $s,s'\in\SCon$ start at $x$ and that the angle measured  them as they leave $x$ is less than $\pi$).
When this is the case, it is easy to find a saddle connection $\widetilde{s}$ with $t(\widetilde{s}) = x$ such that both $\widetilde{s} \rightsquigarrow s$ and $\widetilde{s} \rightsquigarrow s'$.
By \eqref{eq:period-class-transition} we then have that
\[
 j(s)=j(\widetilde{s})+1=j(s') \pmod d.
\]
Now take arbitrary $s,s'$ starting at $x$ and consider their outgoing directions among the finitely many triangulation-edge directions at
$x$. Walking around the circle from $s$ to $s'$ gives a finite sequence of outgoing saddle connections for which their successive angle gaps are strictly less than $\pi$. By the preceding paragraph, $j$ is constant along this sequence. This proves the claim and we may therefore write $j(s)=a_x\in\mathbb Z/d\mathbb Z$ whenever $i(s)=x$.

We now note that for any saddle connection $s$ from $x_1 \in \mathcal X$ to $x_2 \in \mathcal X$. Equation \eqref{eq:period-class-transition} then gives
\[
 a_{x_2}=j(x_2)=j(x_2)+1=a_{x_1}+1 \pmod d.
\]
Applying the same argument to the reversed saddle connection from
$x_2$ to $x_1$ shows that $a_{x_1}=a_{x_2}+1\pmod d$. Hence $d$ divides $2$.

Finally, consider any triangle in the triangulation.
If two of its vertices represent the same singularity $x$, then one
side is a saddle connection from $x$ to itself and $d=1$. Otherwise its vertices
are three distinct singularities $x_1,x_2,x_3 \in \mathcal X$. Applying the relation to
its three oriented sides $x_1\to x_2$, $x_2\to x_3$ and $x_3\to x_1$ implies that $d$ divides $3$.  Since $d$ divides $2$ and $3$ we conclude that $d=1$ and the proof is finished.
\end{proof}

The following result will be key to establish the required non-arithmeticity condition for the local limit theorem.

\begin{lemma} \label{lem:NA}
There does not exist $\tau \in \R \backslash \{0\}$ and $\beta\in\R$ such that
\[
\tau\ell(\gamma)-\beta m(\gamma)\in 2\pi\mathbb Z.
\]
for all closed geodesics $\gamma$.
\end{lemma}
\begin{proof}
If the conclusion failed then the length of every closed geodesic would belong to
$\tau^{-1}(\beta\mathbb Z+2\pi\mathbb Z)$, a rational vector space generated by two real numbers. This contradicts Proposition \ref{prop:geo} (iii). 
\end{proof}

\subsection{Functional analysis}\label{sec.fun}
As discussed in the introduction, we will define a natural operator on $\ell^2(\SCon)$ and study its spectral properties. Here we briefly recall some preliminary material. See Section 7.1 of \cite{EisnerFarkasHaaseNagel2015} for more details.

The space $\ell^2(\SCon)$ is a Hilbert space when equipped with the standard inner product. 
It can also be equipped with a Banach lattice structure:

For real-valued sequences $u,v \in \ell^2(\SCon)$, we write $u \leq v$ if $u(s) \leq v(s)$ for every $s \in \SCon$ (here $u(s),v(s)$ are the entries of $u,v$ at index $s$). For any such $u$ and $v$, their supremum $u \vee v$
and infimum $u \wedge v$ are defined coordinate-wise by
\[
(u \vee v)(s) = \max\{u(s),v(s)\} \ \ \text{ and } \ \  (u \wedge v)(s) = \min\{u(s),v(s)\},
\]
for every $s \in \SCon$. Both sequences
belong to $\ell^2(\SCon)$. The lattice modulus is given by
$|u|(s)=|u(s)|$, which also extends to complex-valued
sequences. This norm is monotone in the sense
that $|u| \leq |v|$ implies $\|u\|_2 \leq \|v\|_2$.
Together with completeness, these properties make
$\ell^2(\SCon)$ a complex Banach lattice
(see \cite[Definition 7.2]{EisnerFarkasHaaseNagel2015}).

A linear subspace $I\subseteq\ell^2(\SCon)$ is called
a lattice ideal if
\[
u \in I,\ \ \  |v|\leq |u| \ \ \ \Longrightarrow \ \ \  v\in I,
\]
see \cite[Definition 7.9]{EisnerFarkasHaaseNagel2015}.
A positive, bounded linear operator $T:\ell^2(\SCon)\to\ell^2(\SCon)$ is said to be
Banach lattice irreducible if its only invariant, closed lattice ideals are
$\{0\}$ and $\ell^2(\SCon)$
\cite[Definition 7.13]{EisnerFarkasHaaseNagel2015}.

Every closed lattice ideal of $\ell^2(\SCon)$ is of the form
\[
I_J=\{u\in\ell^2(\SCon): u(s)=0\text{ for every } s\notin J\}
\]
for some $J\subseteq\SCon.$
Indeed, if $I$ is a closed lattice ideal and $J$ is the set of coordinates at which some element of $I$ is nonzero, then the ideal property implies that $I$ contains every standard basis vector
$e_s$ with $s \in J$. Since $I$ is closed it therefore contains its closed linear span, which is $I_J$.
The reverse inclusion follows from the definition of $J$, and conversely every $I_J$ is a closed lattice ideal.

We are now ready to proceed with the proof of our local limit theorem.

\section{A trace-class operator and its spectral data}\label{sec:spectral}
The proof is based around understanding the properties of the following family of linear operators. Let $(u_s)_{s\in\SCon}$ be the standard orthonormal basis of $\ell^2(\SCon)$. 
\begin{definition}
For $r \in\C$ with $\textnormal{Re}(r)>0$, define the following operator on $\ell^2(\SCon)$: for $s,s' \in \SCon$ set
\[
A_r(s,s')=\id_{\{s\rightsquigarrow s'\}} \exp\!\left(-\frac{r}{2} \left(\ell(s)+\ell(s')\right)\right)
\]
where $\id_{\{s\rightsquigarrow s'\}} = 1$ if $s \rightsquigarrow s'$ and $0$ otherwise. Then for $v \in \ell(\SCon)$ we set
\begin{equation}\label{eq:matrixseries}
 A_r (v) =\sum_{s,s'\in\SCon}   A_r(s,s') \, E_{ss'} (v)
\end{equation}
where $E_{ss'}(v) = \langle v, u_{s'} \rangle u_s$.
\end{definition}
This operator can be thought of as an infinite matrix. The operator is defined precisely so that the  product of the weights around a closed path $w$ equals exactly $e^{-r\ell(w)}$. The definition of this operator is inspired by the operator considered in Section 3 of \cite{CP25}. We now turn our attention to understanding the spectral properties of this result and how it can be used to encode closed geodesics.

For the following result, recall that the trace norm $\|A\|_1$ of a compact operator $A : H \to H$ (from a Hilbert space $H$ to itself) is the sum of the singular values of $A$.

\begin{lemma} \label{lem:trace}
For every $r\in\mathbb C$ with $\textnormal{Re}(r)>0$, the infinite matrix
$A_r$ defines a trace-class operator on
$\ell^2(\SCon)$. More precisely, writing $\|\cdot\|_1$ for the
trace norm and $Z(t)=\sum_{s\in\SCon}e^{-t\ell(s)}$ for $t>0$, we have that
\begin{equation}\label{eq:tracenorm}
 \|A_r\|_1\le Z(\textnormal{Re}(r)/2)^2<\infty.
\end{equation}
The map $r\mapsto A_r$ is holomorphic with values in the space of  trace-class operators on
 $\ell^2(\SCon)$ throughout $\textnormal{Re}(r)>0$.
Moreover, for every integer $n\ge1$,
\begin{equation}\label{eq:traceidentity}
\tr A_r^n =\sum_{\substack{w\text{ a  closed path}\\m(w)=n}}e^{-r\ell(w)},
\end{equation}
and the sum is absolutely convergent.
\end{lemma}
\begin{proof}
Let $F_R=\{s\in\SCon:\ell(s)\le R\}$ for $R>0$. By Proposition \ref{prop:geo}(i), $F_R$ is finite for every $R>0$, and
\begin{equation*}\label{eq:Zbound}
 Z(t)=\sum_{s\in\SCon}e^{-t\ell(s)}<\infty
\end{equation*}
for $t>0$.

Now note that the operators $E_{ss'}$ in \eqref{eq:matrixseries} have exactly one nonzero singular value, equal to $1$. Hence $\|E_{ss'}\|_1=1.$ 
Writing $\alpha=\textnormal{Re}(r)>0$, we then have
\begin{align*}
 \sum_{s,s'\in\SCon} \left\|A_r(s,s') \, E_{ss'} (v)\right\|_1  &\le\sum_{s,s'\in\SCon} e^{-\alpha\ell(s)/2}e^{-\alpha\ell(s')/2} \\
 & =\left(\sum_{s\in\SCon}e^{-\alpha\ell(s)/2}\right)^2 =Z(\alpha/2)^2<\infty.
\end{align*}

For $N\in\mathbb N$, set
\[
 A_r^{(N)}=\sum_{s,s'\in F_N} A_r(s,s') \, E_{ss'} .
\]
Note that $A_r^{(N)}$ is finite rank and hence trace class.
For $M>N$, the triangle inequality gives
\begin{align*}
 \|A_r^{(M)}-A_r^{(N)}\|_1 &\le \sum_{(s,s')\in F_M^2\backslash F_N^2} e^{-\alpha(\ell(s)+\ell(s'))/2} \\
 &\le \sum_{(s,s')\notin F_N^2}
      e^{-\alpha(\ell(s)+\ell(s'))/2} =Z(\alpha/2)^2-  \left(\sum_{s \in F_N}e^{-\alpha\ell(s)/2}\right)^2 \to 0
 \label{eq:sequencetail}
\end{align*}
where the last limit follows from the convergence of $Z(\alpha/2)$ above.
Therefore $(A_r^{(N)})_{N\ge1}$ is a Cauchy sequence in the trace norm.
We then use that $(\mathfrak S_1(\ell^2(\SCon)),\|\cdot\|_1)$, the space of trace class operators on $\ell^2(\SCon)$, is a Banach space and is hence  complete in the trace norm; see
Peterson \cite[Proposition 2.1.6, p. 21]{Peterson13}, as well as
Simon \cite[Chapters 2-3]{Simon05}.  Therefore $A_r^{(N)}$ converges to a trace-class operator. It is easy
to check that $A_r^{(N)}$ tends to $A_r$ as $N\to\infty$. Hence each $A_r$ is trace class as required.\\

Analyticity is easy to check from the above. Indeed,  fix $r_0$ with $\alpha_0=\textnormal{Re}(r_0)>0$, and work in the disk
$|r-r_0|<\alpha_0/2$. There $\textnormal{Re}(r)>\alpha_0/2$, and the trace-norm sum of the terms in \eqref{eq:matrixseries} is bounded above by $Z(\alpha_0/4)^2$. Every individual operator that is summed in  \eqref{eq:matrixseries} is entire in $r$.
The Weierstrass theorem for Banach-space-valued holomorphic functions (Theorem 3.1 in \cite{weierstrass})
therefore proves local uniform convergence and analyticity in $\mathfrak S_1(\ell^2(\SCon))$. \\

We now verify the trace formula.
Note that for any fixed $n\ge1$ and $j \ge 1$,
\begin{equation}\label{eq:tracecycles}
 \tr\left((A_r^{(j)})^n\right) =\sum_{s_1,\ldots,s_n\in F_j}  \prod_{k=1}^n\id_{\{s_k\rightsquigarrow s_{k+1}\}} \exp\!\left(-\frac{r}{2} \sum_{k=1}^n(\ell(s_k)+\ell(s_{k+1}))\right)  =\sum_{\substack{(s_1,\ldots,s_n)\in F_j^n\\ s_k\rightsquigarrow s_{k+1}\ (\mathrm{closed})}}   e^{-r\sum_{k=1}^n\ell(s_k)},  
\end{equation}
where $s_{n+1}=s_1$. The  second equality follows because
each $\ell(s_k)$ occurs in two neighbouring factors with coefficient $-r/2$. 

The right-hand sums converge absolutely as $j\to\infty$, since
\begin{equation}\label{eq:cyclesabsolute}
 \sum_{\substack{s_1,\ldots, s_n\in\SCon\\
                 s_k\rightsquigarrow s_{k+1}\ (\mathrm{closed})}} \left|e^{-r\sum_{k=1}^n\ell(s_k)}\right|  \le\sum_{s_1,\ldots,s_n\in\SCon} \prod_{k=1}^ne^{-\alpha\ell(s_k)}=Z(\alpha)^n<\infty.
\end{equation}

It remains to justify passing to the limit on the left hand  side of \eqref{eq:tracecycles}. Put $B_N=A_r^{(N)}$ and
$B=A_r$ so that $\|B_N-B\|_1\to0$ as $N\to\infty$.
For every fixed positive integer $n$, we will use the following elementary algebraic identity:
\begin{equation}\label{eq:powerstelescope}
 B_N^n-B^n
 =\sum_{j=0}^{n-1}B_N^{\,n-1-j}(B_N-B)B^j.
\end{equation}
Recall that the trace-class operators satisfy the following:
if $T$ is trace class and $S,R$ are bounded, then
\begin{equation}\label{eq:twosided}
 \|STR\|_1\le\|S\|_{\mathrm{op}}\,\|T\|_1\,   \|R\|_{\mathrm{op}},
\end{equation}
see Simon \cite[Chapters 2-3]{Simon05}. Here $\|\cdot\|_{\mathrm{op}}$ is the operator norm. Applying
\eqref{eq:twosided} to each term of \eqref{eq:powerstelescope}
gives
\begin{align*}
 \|B_N^n-B^n\|_1
 &\le\sum_{j=0}^{n-1}
   \|B_N\|_{\mathrm{op}}^{\,n-1-j}
   \|B_N-B\|_1
   \|B\|_{\mathrm{op}}^{\,j} \\
 &\le nM^{n-1}\|B_N-B\|_1,
\end{align*}
where $M=Z(\alpha/2)^2$. The last inequality uses
$\|B_N\|_{\mathrm{op}}\le\|B_N\|_1\le M$ and
$\|B\|_{\mathrm{op}}\le\|B\|_1\le M$,
which follow from \eqref{eq:tracenorm}.
Because $M$ and $n$ are fixed, this proves that $B_N^n\to B^n$ in trace norm.

For a trace-class operator $T$, the trace is a continuous linear
functional satisfying $|\tr T|\le\|T\|_1$; see
Simon \cite[Chapters 2-3]{Simon05}. Hence
\begin{equation}\label{eq:tracepowerlimit}
 \left|\tr((A_r^{(N)})^n)-\tr(A_r^n)\right|  \le\left\|(A_r^{(N)})^n-A_r^n\right\|_1 \longrightarrow0.
\end{equation}
On the other hand, each fixed closed path 
$(s_1,\ldots,s_n)$ lies in $F_N^n$ for all sufficiently large $N$.
Since \eqref{eq:cyclesabsolute} dominates the absolute values
of all the terms by a summable function on $\SCon^n$,
the dominated convergence theorem gives that
\[
 \lim_{N\to\infty}  \sum_{\substack{(s_1,\ldots,s_n)\in F_N^n\\
                   s_k\rightsquigarrow s_{k+1}\ \mathrm{(closed)}}}
 e^{-r\sum_k\ell(s_k)}
 =\sum_{\substack{(s_1,\ldots,s_n)\in\SCon^n\\
                   s_k\rightsquigarrow s_{k+1}\ \mathrm{(closed)}}} e^{-r\sum_k\ell(s_k)}.
\]
Taking limits in \eqref{eq:tracecycles}, using \eqref{eq:tracepowerlimit} for its left side and  the previous display equation for its right side, proves \eqref{eq:traceidentity}.
\end{proof}

We now show that the operators $A_r$ are Banach lattice irreducible as defined in Section \ref{sec.fun}. For $J \subset \SCon$ let $I_J$ denote the ideal defined there.

\begin{proposition}\label{prop.bli}
For every real $r>0$, the operator $A_r:\ell^2(\SCon)\to\ell^2(\SCon)$ is Banach lattice irreducible.
\end{proposition}

\begin{proof}
Since $r>0$, all matrix coefficients of $A_r$ are
nonnegative, so $A_r$ is a positive operator.
Suppose, for a contradiction, that there exists a
nontrivial closed lattice ideal $I_J$ invariant under
$A_r$ (for some $J \subset \SCon$). Choose $s'\in J$ and $s\in\SCon\backslash J$.
By Proposition 2.3(ii), there is a directed path $s=s_0\rightsquigarrow s_1\rightsquigarrow \cdots\rightsquigarrow s_n= s'.$
Since each transition along this path has strictly positive weight, we obtain
\[
(A_r^n)(s,s') \geq \prod_{k=0}^{n-1} A_r(s_k,s_{k+1})>0.
\]
Let $u_{s'}$ denote the standard coordinate vector
corresponding to $s'$. Since $ s'\in J$, we have
$u_{s'}\in I_J$. By the invariance of $I_J$, $A_r^n u_{s'}\in I_J.$
Then, as $s \notin J$, $(A_r^n u_{s'})(s)=0.$
On the other hand,
\[
(A_r^n u_{s'})(s)=(A_r^n)(s,s')>0,
\]
which is a contradiction. Therefore no nontrivial closed lattice ideal is invariant under $A_r$, and hence $A_r$ is Banach lattice irreducible.
\end{proof}

We can now deduce good properties of the spectrum of $A_r$ using the aforementioned Perron-Frobenius theorem \cite[Theorem 1.5]{CWW20}.

\begin{proposition}\label{prop:PF}
There is a unique $h>0$ such that $\rho(A_h)=1$. On a real neighbourhood $I$ of $h$, $A_r$ has a simple positive spectral radius $\lambda(r)$, depending real-analytically on $r$, and
\begin{equation}\label{eq:spectrumper}
\textnormal{spec}(A_r)\cap\{z:|z|=\lambda(r)\}
=\{\lambda(r)\},
\end{equation}
is simple. Lastly,
\begin{equation*}
 a=-\left.\frac{d}{ds}\log\lambda(s)\right|_{s=h}  \ \ \text{ and } \ \  v=\left.\frac{d^2}{ds^2}\log\lambda(s)\right|_{s=h} 
\end{equation*}
are both strictly positive.
\end{proposition}

\begin{proof}
\textit{Perron-Frobenius theory and the spectral gap.}
By Lemma \ref{lem:trace} and Proposition \ref{prop.bli}, $A_r$ is a
positive, compact and Banach lattice irreducible operator. Therefore by \cite[Theorem 1.5]{CWW20}, its spectral
radius $\lambda(r): = \rho(A_r) >0$ is an algebraically simple
eigenvalue, with strictly positive right and left
eigenvectors $u_r,v_r\in\ell^2(\SCon)$:
\begin{equation*}
  A_r u_r=\lambda(r)u_r, \ \  A_r^*v_r=\lambda(r)v_r, \ \ u_r(s),v_r(s)>0
\end{equation*}
for $s \in \SCon$. By compactness, the spectrum of $A_r$, other than those with modulus $\lambda(r)$ is contained in a disk of radius strictly less that $\lambda(r)$.

We now verify that $\mathscr G_X$ having period one excludes other
eigenvalues of modulus $\lambda(r)$.. Suppose
$A_rz=\mu z$ with $z\ne0$ and $|\mu|=\lambda(r)$.
The coordinate-wise triangle inequality gives
\[
  \lambda(r)|z|=|A_rz|\le A_r|z|.
\]
We then take the product of this inequality with the strictly positive vector $v_r$.
Using $A_r^*v_r=\lambda(r)v_r$, we see that the nonnegative
vector $A_r|z|-\lambda(r)|z|$ has zero inner product with $v_r$.
It must therefore vanish coordinate-wise. By simplicity of the eigenvalue, $|z|=c\,u_r$ for some $c>0$.
After rescaling, write
$z(s)=u_r(s)\zeta(s)$ with $|\zeta(s)|=1$.
With $\eta:=\mu/\lambda(r)$, the eigenvector equations imply
\[
  \sum_{s' \in \SCon} A_r(s,s')u_r(s')\zeta(s') =\eta\zeta(s)\sum_{s' \in \SCon} A_r(s,s')u_r(s').
\]
The coefficients on the right-hand sum are positive, and the
weighted average of the unit complex numbers $\zeta(s')$ again has
modulus one. This forces
\begin{equation*}
  \zeta(s')=\eta\zeta(s)
   \ \ \ \text{for every }s\rightsquigarrow s'.
\end{equation*}
Going around any closed path of saddle connection count $m(w)$ gives
$\eta^{m(w)}=1$. The gcd of all these counts is one by
Proposition \ref{prop:aperiodic}, hence $\eta=1$.
It follows that
\[
  \textnormal{spec}(A_r)\cap\{|z|=\lambda(r)\}=\{\lambda(r)\}.
\]
Since $A_r$ is compact and $\lambda(r)$ is a simple isolated
eigenvalue, analytic perturbation theory \cite[Chapter VII, §1.3, Theorem 1.8]{Kato} shows that
$r\mapsto\lambda(r)$ is real-analytic on $(0,\infty)$.

\medskip
\textit{Existence and uniqueness of $h$:}
For $\alpha>0$ we have the entry-wise inequality
\[
  0\le A_{r+\alpha}\le e^{-\alpha\ell_*}A_r
\]
therefore
\[
  \lambda(r+\alpha)\le e^{-\alpha\ell_*}\lambda(r)<\lambda(r).
\]
Thus $\lambda$ is strictly decreasing. Lemma \ref{lem:trace} also gives
\[
  \lambda(r)\le\|A_r\|_1\le Z(r/2)^2 \to 0
\]
as $r\to\infty$ where the convergence follows by dominated convergence, using
$Z(t_0)<\infty$ for any fixed $t_0>0$.

It is easy to check that $\lambda(r)>1$ for sufficiently small $r>0$.
By continuity, strict monotonicity, and $\lambda(r)\to0$ at infinity, there is a unique $h>0$
for which $\lambda(h)=1$.\\

To prove positivity of the constants we apply the work in Cantrell-Pollicott. See Appendix \ref{sec.ident} that identifies $a$ and $v$ with the mean and variance corresponding to the counting central limit theorem: Theorem 1.2 in \cite{CP25}.
\end{proof}

\begin{definition}
We set 
\[
 \Lambda=\frac{1}{a} \ \ \ \text{ and } \ \ \  \sigma^2=\frac{v}{a^3}.
\]
\end{definition}

\begin{remark}
These will turn out to be the mean and variance of our local limit theorem and they coincide with the mean and variance in Theorem \ref{thm.sing} above. We also note that it is simple to check that $h>0$ in the above proposition is the exponential growth rate of $\#\Prim(T)$ as $T\to\infty$.
\end{remark}

The following result controls the spectral radius of the operators $A_r$ away from $r = h$.
\begin{lemma}\label{lem:frequency}
For all real $r\in I$ and $\tau\ne0$,
\begin{equation}\label{eq:strict}
 \rho(A_{r+i\tau})<\lambda(r).
\end{equation}
For any compact $K\subset\R\backslash\{0\}$ there exist $C_K>0$ and $0<q_K<1$ such that
\[
\sup_{r\in I',\,\tau\in K} \left\|\lambda(r)^{-n}A_{r+i\tau}^n\right\|_1 \le C_Kq_K^n
\]
for all $n\ge1$, where $I'\subset I$ is any sufficiently small closed interval containing $h$.
\end{lemma}
\begin{proof}
Fix $r\in I$ and $\tau\ne0$, and abbreviate
\[
 B=A_{r+i\tau}, \ \ \  A=A_r, \ \ \  \lambda=\lambda(r), \ \ \  g(s,s')=\tfrac{1}{2}\left(\ell(s)+\ell(s')\right).
\]
Then $B(s,s')=A(s,s')e^{-i\tau g(s,s')}$ and $|B(s,s')|=A(s,s')$ for all $s,s'$.  It follows that
for $z\in\ell^2(\SCon)$ and every $n\ge1$ we have that $ |B^nz|\le A^n|z|$ coordinate wise and so
\[
 \|B^n\|_{\mathrm{op}}\le\|A^n\|_{\mathrm{op}}.
\]
The spectral-radius formula therefore gives $\rho(B)\le \rho(A)=\lambda$.

We now want to understand the equality case.  This is essentially the countable-state analogue of
Wielandt's theorem for an irreducible nonnegative matrix
\cite[Chapter 8]{HornJohnson13}.  Since $B$ is trace class by Lemma \ref{lem:trace}, it is compact. Now,  if $\rho(B)=\lambda>0$,
then there is an eigenvector $z \ne 0$ and some $\beta\in \mathbb{R}$ such that
\begin{equation}\label{eq:eigeneq}
 Bz=\lambda e^{-i\beta}z.
\end{equation}
Let $u_r,v_r\in\ell^2(\SCon)$ be the strictly positive right and left
lead eigenvectors of $A$, so that
$Au_r=\lambda u_r$ and $A^*v_r=\lambda v_r$ (Proposition \ref{prop:PF}).
Taking absolute values of \eqref{eq:eigeneq} gives
\[
 \lambda|z|=|Bz|\le A|z|.
\]
The difference $A|z|-\lambda|z|$ belongs to $\ell^2$, is
coordinate-wise nonnegative, and has zero product with $v_r$:
\[
 \left\langle v_r,A|z|-\lambda|z|\right\rangle =\left\langle A^*v_r-\lambda v_r,|z|\right\rangle=0.
\]
Because $v_r(s)>0$ for every state $s$, each coordinate of this
difference must vanish.  Thherefore $A|z|=\lambda|z|$ and simplicity of the
lead eigenvalue implies $|z|=cu_r$ for some $c>0$.  Rescale
$z$ so that $c=1$, and write $z_s=u_r(s)\zeta_s$, where
$|\zeta_s|=1$ for every $s$.

The $s$th coordinate of \eqref{eq:eigeneq} is now
\begin{equation}\label{eq:row}
 \sum_{s':s\rightsquigarrow s'} A(s,s')u_r(s')e^{-i\tau g(s,s')}\zeta_{s'} =\lambda u_r(s)e^{-i\beta}\zeta_s,
\end{equation}
whereas the eigenvalue equation yields
$\sum_{s':s\rightsquigarrow s'}A(s,s')u_r(s')=\lambda u_r(s)$.
Now, all of the terms in the sum in left hand side of \eqref{eq:row} must have the same argument as the right-hand
side.  Since $A(s,s')u_r(s')>0$, on every admissible transition $s\rightsquigarrow s'$, it follows
that
\begin{equation}\label{eq:coboundary}
 e^{-i\tau g(s,s')}\zeta_{s'}=e^{-i\beta}\zeta_s.
\end{equation}
For a closed path $w=(s_1,\ldots,s_m)$, multiply \eqref{eq:coboundary}
along its $m$ transitions, with $s_{m+1}=s_1$.  The $\zeta_s$ factors cancel, leaving
\[
e^{-i\tau\ell(w)}=e^{-i\beta m(w)},
\]
or equivalently $\tau\ell(w)-\beta m(w)\in2\pi\mathbb Z$ for every closed path $w$.
This contradicts Lemma \ref{lem:NA}. Therefore $\rho(B)<\lambda$, proving \eqref{eq:strict}.\\

For the final statement, let $\mathcal K=I'\times K$ and set
$D_{r,\tau}=\lambda(r)^{-1}A_{r+i\tau}$.
The conclusion follows from the facts that  the family $D_{r,\tau}$ depends continuously (in fact analytically by Lemma \eqref{lem:trace}) on $(r,\tau)$ in trace
norm on the compact set $\mathcal K$, and every member has spectral radius less than one. 
\end{proof}

\section{Deducing the local limit theorem}\label{sec:fixed}
We are now ready to deduce the local limit theorem. We follow the arguments in \cite{Lalley1987}. For $n\ge1$ and $T>0$, let $A_n(T)$ count closed saddle paths:
\begin{equation*}
A_n(T)=\#\{w:\ w\text{ a  closed saddle path},\ m(w)=n,\ \ell(w)<T\}.
\end{equation*}
By Proposition \ref{prop:aperiodic}, $A_n(T)$ is positive for every $n$ and sufficiently large $T$ depending on $n$.

\begin{lemma}\label{lem:stone}
Let $I'\subset I$ be a compact interval containing $h$. For $r\in I'$ put $a_r=-(\log\lambda)'(r)$ and $v_r=(\log\lambda)''(r)>0$. Uniformly for $r\in I'$, $T\in\R$,
\begin{equation}\label{eq:stone}
A_n(T)=\frac{e^{rT}\lambda(r)^n}{r\sqrt{2\pi v_r n}} \exp\!\left(-\frac{(T-a_rn)^2}{2v_rn}\right) +o\, \!\left(\frac{e^{rT}\lambda(r)^n}{\sqrt n}\right)
\end{equation}
as $n\to\infty$.
In particular, there are constants $C>0$ and $n_0\in\N$ independent of $T,r$ such that
\begin{equation}\label{eq:uniformupper}
A_n(T)\leq C\,\frac{e^{rT}\lambda(r)^n}{\sqrt n} \ \ \text{ for }  n\ge n_0,\ r\in I'.
\end{equation}
\end{lemma}

\begin{proof}

Fix a compact interval $I'\subset I$ containing $h$.  Throughout the
proof, all constants and error terms are uniform for $r\in I'$.
Set, for $r \in I'$, $T>0$ and $n\ge 1$,
\[
 L(r)=\log\lambda(r), \ \ a_r=- L'(r),\  v_r=L''(r) \ \text{ and } \ \  G_{n,r}(T)=\frac{1}{\sqrt{2\pi v_r n}}  \exp\!\left(-\frac{(T-a_r n)^2}{2v_r n}\right).
\]
By Proposition \ref{prop:PF} and compactness of $I'$, there are constants
$r_->0$ and $v_->0$ with $r\ge r_-$ and $v_r\ge v_-$ for $r\in I'$.
In particular $0\le G_{n,r}(T)\le(2\pi v_-n)^{-1/2}$.

Let $\nu_{n,r}$ be the  measure that puts mass $\lambda(r)^{-n}e^{-r\ell(w)}$ at
$\ell(w)$ for each closed path $w$ with $m(w)=n$.
Lemma \ref{lem:trace} gives the identity
\begin{equation*}
 \widehat\nu_{n,r}(\tau)
 :=\int_{\mathbb R}e^{i\tau t}\,\nu_{n,r}(dt)
 =\lambda(r)^{-n} \, \tr A_{r-i\tau}^{\,n}
\end{equation*}
for the Fourier transform $\widehat{\nu}_{n,r}$.

Now, by Proposition \ref{prop:PF}, the simple leading eigenvalue admits an analytic
continuation $\lambda(z)$ to a complex neighbourhood of $I'$.
After reducing $\epsilon>0$ if necessary, there are analytic spectral projections $\Pi_{r,\tau}$ and remainders $N_{r,\tau}$ such
that, for $r\in I'$ and $|\tau|\le\epsilon$,
\[
 A_{r-i\tau}=\lambda(r-i\tau)\Pi_{r,\tau}+N_{r,\tau}, \ \ \Pi_{r,\tau}N_{r,\tau}=N_{r,\tau}\Pi_{r,\tau}=0
\]
see \cite[Theorem 1.7, pp. 368–369]{Kato}.
The remainder spectra lie, uniformly, in a disk of radius
$q_0\lambda(r)$, with some $q_0<1$.  Since the operators depend
continuously in trace norm we have that
\begin{equation}\label{eq:lnear}
 \widehat\nu_{n,r}(\tau) =\left(\frac{\lambda(r-i\tau)}{\lambda(r)}\right)^n  +O(q^n)
\end{equation}
uniformly in $ |\tau|\le\epsilon.$
Using Taylor's theorem, uniformly in $r$, gives
\begin{equation}\label{eq:taylor}
 \log\frac{\lambda(r-i\tau)}{\lambda(r)}
 = i a_r\tau-\frac12 v_r\tau^2+O(|\tau|^3)
\end{equation}
and because $v_r\ge v_->0$, we can reduce $\epsilon$ further so that
\begin{equation}\label{eq:gaussian-decay}
 \textnormal{Re}\log\frac{\lambda(r-i\tau)}{\lambda(r)}
 \le-c\tau^2
\end{equation}
for all $|\tau|\le\epsilon$ and  $c>0$ independent of $\tau$.

For every $0<\epsilon\le B$, Lemma \ref{lem:frequency} provides constants
$C_B>0$ and $0<q_B<1$ such that
\begin{equation}\label{eq:away}
 \sup_{\substack{r\in I'\\\epsilon\le|\tau|\le B}}
 |\widehat\nu_{n,r}(\tau)|\le C_Bq_B^n.
\end{equation}

Let $f$ be a Schwartz function whose Fourier transform
$\widehat f(\tau)=\int e^{-i\tau u}f(u)\,du$ has compact support.
Fourier inversion  gives
\begin{equation}\label{eq:invert}
 \int f(t-T)\,\nu_{n,r}(dt) =\frac1{2\pi}\int_{\mathbb R} \widehat f(\tau)e^{-i\tau T}\widehat\nu_{n,r}(\tau)\,d\tau.
\end{equation}

The portion of this integral with $|\tau|\ge\epsilon$ is $O(q^n)$ for some $0 < q < 1$ by \eqref{eq:away}. Also, when $|\tau|\le\epsilon$, \eqref{eq:lnear} - \eqref{eq:gaussian-decay} and the substitution $\tau=u/\sqrt n$ give, uniformly for $r\in I'$,
\[
 \int_{-\epsilon}^{\epsilon} \left|\left(\frac{\lambda(r-i\tau)}{\lambda(r)}\right)^n -e^{ia_rn\tau-v_rn\tau^2/2}\right|d\tau = o(n^{-1/2})
\]
as $n\to\infty$.
Indeed, after the substitution, the integrals converge to zero uniformly in $r$ for each fixed $u$, and each are bounded above by $2e^{-c'u^2}$ for some $c' >0$.  Also, using continuity of $\widehat f$ at $0$, we can replace $\widehat f(\tau)$ by $\widehat f(0)=\int f$ at the cost of introducing an error that it at most $o(n^{-1/2})$.

The factor $e^{-i\tau T}$ has absolute value one, so every error is uniform in $T\in\mathbb R$.  Substituting the above into \eqref{eq:invert} gives \begin{equation}\label{eq:band}
 \sup_{\substack{r\in I'\\T\in\mathbb R}} \sqrt n\left| \int f(t-T)\,\nu_{n,r}(dt) -G_{n,r}(T)\int_{\mathbb R}f(u)\,du\right|
\to 0.
\end{equation}

A standard smoothing step (which we omit) in the proof of local limit theorems (see \cite[Section 7]{Lalley1987} for example) extends this statement (first to indicator functions and then) to the Riemann-integrable function $f_r(u)=e^{ru} \, \id_{\{u<0\}}$, uniformly for $r\in I'$.  With this choice, we have the exact identity
\[
 \int f_r(t-T)\,\nu_{n,r}(dt) =e^{-rT}\lambda(r)^{-n}A_n(T).
\]
Applying this extension of \eqref{eq:band} and multiplying by
$e^{rT}\lambda(r)^n$ proves
\[
 A_n(T)
 =\frac{e^{rT}\lambda(r)^n}{r\sqrt{2\pi v_r n}}
   \exp\!\left(-\frac{(T-a_rn)^2}{2v_rn}\right)
   +o\!\left(\frac{e^{rT}\lambda(r)^n}{\sqrt n}\right),
\]
as $n\to\infty$ where the little-$o$ is uniform in $r\in I'$ and $T\in\mathbb R$. This also gives \eqref{eq:uniformupper} for all sufficiently large $n$.
\end{proof}

\begin{lemma}\label{lem:tails}
There exist $\delta,c,C>0$ such that, whenever $n\in\mathbb N$ and $|n-T/a|\leq\delta T$,
\begin{equation}\label{eq:gausstail}
A_n(T)\leq\frac{Ce^{hT}}{\sqrt T} \exp\!\left(-c \, \frac{(n-T/a)^2}{T}\right).
\end{equation}
Furthermore, for some $c_\delta>0$,
\begin{equation}\label{eq:largetail}
\sum_{n\ge1:\ |n-T/a|>\delta T}A_n(T)=O(e^{(h-c_\delta)T})
\end{equation}
as $T\to\infty$.
\end{lemma}
\begin{proof}
As before, write $L(r)=\log\lambda(r), a_r=-L'(r)$ and $v_r=L''(r)$ so that  $L(h)=0$, $a_h=a>0$, and $v_h=v>0$.
Shrinking the interval $I$ from Proposition \ref{prop:PF} if necessary,
choose a compact interval $J\subset I$, with $h$ in its interior,
and constants $0<v_-\le v_+$ such that
\begin{equation}\label{eq:v-bounds}
 v_-\le L''(r)\le v_+
\end{equation}
for $r \in J$. In particular $r\mapsto a_r$ is strictly decreasing on $J$.\\

Choose $\delta>0$ sufficiently small that whenever
\[
\left|n-\frac{T}{a}\right|\le\delta T,
\]
there are constants $c_0,C_0>0$ with
\begin{equation}\label{eq:comparable}
 c_0T\le n\le C_0T
\end{equation}
and furthermore  $T/n$ belongs to the interior of $\{a_r : r \in J\}$. When this is the case there is a unique $r=r(n,T)\in J$
such that $ a_r= T/n$ or equivalently $L'(r)=-T/n$.
By the uniform bound \eqref{eq:uniformupper} in Lemma \ref{lem:stone},
\begin{equation}\label{eq:tilted-count}
 A_n(T)\le \frac{C}{\sqrt n}\exp\left(rT+nL(r)\right).
\end{equation}
We now estimate the exponent $rT+nL(r)$. Since $T=na_r=-nL'(r)$,
\begin{align}
 rT+nL(r)-hT &=n\left(L(r)-(r-h)L'(r)-L(h)\right)\notag\\
 &=-n\int_h^r(t-h)L''(t)\,dt\notag\\
 &\le -\frac{v_-n}{2}(r-h)^2.\label{eq:exponent}
\end{align}
The integral is nonnegative for either $r>h$ or $r<h$.
Also, the mean-value theorem and \eqref{eq:v-bounds} give
\[
 \left|\frac{T}{n}-a\right| =|a_r-a|\le v_+|r-h|,
\]
and therefore
\begin{equation}\label{eq:exponent2}
 rT+nL(r)
 \le hT-\frac{v_-}{2v_+^2}\frac{(T-an)^2}{n}.
\end{equation}
Now $T-an=-a(n-T/a)$ and using \eqref{eq:comparable} in \eqref{eq:tilted-count}-\eqref{eq:exponent2},
we obtain constants $C, c >0$ such that
\[
 A_n(T)\le\frac{Ce^{hT}}{\sqrt T}
 \exp\!\left(-c\frac{(n-T/a)^2}{T}\right)
\]
whenever $|n-T/a|\le\delta T$ and $T$ is sufficiently large.
This proves \eqref{eq:gausstail}.\\

We now show the furthermore statement.
First we record an estimate that is uniform over all
$n\ge1$ for each fixed real $r>0$:
\begin{equation}\label{eq:trace-powers}
 \tr A_r^n\le C_r\lambda(r)^n.
\end{equation}
This follows as in the proof of Lemma \ref{lem:stone}.

Choose $\varepsilon>0$ small enough that
$r_-=h-\varepsilon>0$ and $r_+=h+\varepsilon$ lie in $J$.
Strict monotonicity of $\lambda$ gives
\[
 L_-=L(r_-)>0, \ \ \  L_+=L(r_+)<0.
\]
By Taylor's formula at $h$,
\begin{equation}\label{eq:fixed-tilts}
 L(h\pm\varepsilon)=\mp a\varepsilon+O(\varepsilon^2).
\end{equation}
We may choose $\varepsilon$ sufficiently small (depending on
$\delta$) so that for some $\kappa>0$,
\begin{equation}\label{eq:low-rate}
 -\varepsilon+(a^{-1}-\delta)L_- \le-\kappa \ \ \text{ and } \ \  \varepsilon+(a^{-1}+\delta)L_+\le-\kappa.  
\end{equation}
Indeed, both left sides equal $-a\delta\varepsilon
+O(\varepsilon^2)$, by \eqref{eq:fixed-tilts}.

For every $r>0$ the trace formula gives the bound
\begin{equation}\label{eq:markov}
 A_n(T)\le e^{rT}\tr A_r^n \le C_r\exp\left(rT+nL(r)\right).
\end{equation}
If $n/T\le a^{-1}-\delta$, let $r=r_-$.
Since $L_->0$, equations
\eqref{eq:low-rate} and \eqref{eq:markov} imply
\[
 A_n(T)\le C_- e^{T(r_-+(n/T)L_-)}  \le C_-e^{(h-\kappa)T}.
\]
If $n/T\ge a^{-1}+\delta$ instead set $r=r_+$.
Since $L_+<0$, the analogous estimate gives
\[
 A_n(T)\le C_+e^{T[r_++(n/T)L_+]} \le C_+e^{(h-\kappa)T}.
\]
Therefore the same exponential bound holds with constants independent of $n,T$ satisfying
\[
\left|n-\frac{T}{a}\right| > \delta T.
\]

Finally, each saddle connection has length at least $\ell_*$, so
$A_n(T)=0$ for $n\ell_*\ge T$. There are at most $\lceil T/\ell_*\rceil$ potentially nonzero terms and so
\[
 \sum_{n\ge1: \, |n-T/a|>\delta T} A_n(T) \le C(1+T)e^{(h-\kappa)T} = O(e^{(h-c_\delta)T})
\]
for any fixed $0<c_\delta<\kappa$.
This proves \eqref{eq:largetail}, concluding the proof.
\end{proof}

\subsection{Concluding the proof}
For a prime orbit consisting of $n$ saddle connections, there are exactly $n$ distinct closed paths obtained by choosing the starting saddle connection. Therefore, writing $B_n(T)$ for the number of prime oriented closed geodesics with $n$ saddle connections, we have
\[
B_n(T)=\frac{A_n(T)-R_n(T)}{n},
\]
where $R_n(T)$ counts non prime closed saddle paths.
 The following simple result shows that non-prime orbits are negligible in the count. 
\begin{lemma}\label{lem:non prime}
There is a constant $C>0$ with
\[
\sum_{n\geq1}R_n(T)\leq CT e^{hT/2}
\]
for all $T >1$.  Here $h>0$ denotes the exponential growth rate of $\#\Prim(T)$.
\end{lemma}
\begin{proof}
Every non-prime closed path is a $k \ge 2$ repetition of a prime closed path of length at most $T/k$.
Therefore
\[
\sum_{n \ge 1} R_n(T)\leq\sum_{2\leq k\leq T/\ell_*}A(T/k).
\]
We then recall that there exists $C >0$ such that $A(t)\leq Ce^{ht}$ for all $t\geq0$. Hence, in the above, there are $O(T)$ terms and each is $O(e^{hT/2})$ as $T\to\infty$. The result follows.
\end{proof}

We can now conclude the proof by combining our above work.

\begin{proof}[Proof of Theorem \ref{thm:main}]
For $n\ge 1$, $T>0$ and $u \in \R$ we set
\[
u_n=\frac{n-\Lambda T}{\sqrt T} =\frac{n- T/a}{\sqrt T} ,\ \ \  \ \sigma^2=\frac{v}{a^3},\ \ \  \text{ and } \ \  \varphi(u)=\frac{1}{\sqrt{2\pi\sigma^2}}  e^{-u^2/(2\sigma^2)}.
\]
Also fix $R>0$. 

When $|u_n|\le R$, we have
$n=T/a+O_R(\sqrt T)$ and so $n$ is  comparable with $T$.
Applying Lemma \ref{lem:stone} with $r=h$ (the exponential growth rate of $\#\Prim(T)$) and using
$L(h)=0$, we find, uniformly for these $n$,
\begin{equation}\label{eq:basic}
 A_n(T)=\frac{e^{hT}}{h\sqrt{2\pi vn}}  \exp\!\left(-\frac{(T-an)^2}{2vn}\right) +o_R\!\left(\frac{e^{hT}}{\sqrt T}\right).
\end{equation}
as $T\to\infty$. This is exploiting the uniformity in $T$ in Lemma \ref{lem:stone}.
Since $T-an=-au_n\sqrt T$ and $n=(T/a)(1+au_n/\sqrt T)$ we have that
\begin{align*}
 \frac{(T-an)^2}{2vn}
  &=\frac{a^3u_n^2}{2v}\left(1+O_R(T^{-1/2})\right) =\frac{u_n^2}{2\sigma^2}+O_R(T^{-1/2}); \ \ \text{ and, } \\  \frac{1}{h\sqrt{2\pi vn}}
  &=\frac{1}{ha\sqrt T}\,\frac{1}{\sqrt{2\pi\sigma^2}}  \left(1+O_R(T^{-1/2})\right).
\end{align*}
Substituting these into \eqref{eq:basic} shows that
\begin{equation}\label{eq:central}
 A_n(T)=\frac{e^{hT}}{ha\sqrt T}\, \varphi(u_n) +o_R\!\left(\frac{e^{hT}}{\sqrt T}\right) \ \ \ \text{ for } \ \ \ |u_n|\le R,
\end{equation}
where the error is uniform over all such integers $n$.

Dividing \eqref{eq:central} by $n$, subtracting the non prime contribution (which are uniformly, exponentially negligible by Lemma \ref{lem:non prime}), and then dividing by $B(T)$ gives
\[
\sqrt T\frac{B_n(T)}{B(T)} =\varphi(u_n) +o(1)
\]
as $T\to\infty$ uniformly for all $|u_n| < R$.\\

To conclude we need to control the $n\ge 1$ with $|u_n| \ge R$.
To that end, let $\delta >0$ be as in Lemma \ref{lem:tails}.
When $R  < |u_n| \le \delta \sqrt T$ we have that 
\[
\sqrt T \, \frac{B_n(T)}{B(T)}\leq Ce^{-cR^2}
\]
for some $C,c >0$ by  \eqref{eq:gausstail}. For $ |u_n| > \delta \sqrt T$ we have an exponentially small upper bound by \eqref{eq:largetail}.\\

Combining all of the above with the fact that
\[
\sup_{|u_n| > R} \exp\!\left(-\frac{(n-\Lambda T)^2}{2\sigma^2T}\right) \le C' e^{-c'R^2}
\]
for some $C',c' >0$, we deduce that 
\[
\limsup_{T\to\infty}\sup_{n\in\N}\left|\sqrt T\, \frac{B_n(T)}{B(T)} -\frac{1}{\sqrt{2\pi\sigma^2}} \exp \, \left(-\frac{(n-\Lambda T)^2}{2\sigma^2T}\right) \right| \le Ce^{-cR^2}  +  C' e^{-c'R^2}.
\]
Letting $R\to\infty$ gives the required asymptotic for $B_n(T)$, concluding the proof.
\end{proof}

\appendix

\section{Identification and positivity of the constants}\label{sec.ident}

In this section we identify the constants $a, v$ from Proposition \ref{prop:PF} with the mean and variance appearing in Theorem \ref{thm.sing}.

Let $\mathcal P_x$ be the collection of nonempty finite saddle-connection
paths starting at a fixed singularity $x$ as defined in the introduction. 
For $p=(s_1,\ldots,s_n)$ write
$\ell(p)=\sum_j\ell(s_j)$ and $m(p)=n$.
By Theorem \ref{thm.sing} there is a central limit theorem with mean $\Lambda$ and $\sigma^2$ which are both strictly positive.

For real $r>0$ and $t$ near zero set $B_{r,t}:=e^t A_r$ where $A_r$ is as in Section \ref{sec:spectral}.
For $r>0, x \in \mathcal X$ and $s \in \SCon$ let
\[
\xi_{r,x}(s):=\id_{\{i(s)=x\}}e^{-r\ell(s)/2},\ \ \ \  \eta_r(s):=e^{-r\ell(s)/2}.
\]
These determine points in $\ell^2(\SCon)$ by Proposition \ref{prop:geo}(i).
By construction we have the exact formal identity
\begin{equation}\label{eq:pathequal}
\sum_{p\in\mathcal P_x}e^{-r\ell(p)+tm(p)} =e^t\sum_{k\geq0} \left\langle \xi_{r,x},B_{r,t}^{\,k}\eta_r\right\rangle,
\end{equation}
where both sides could be infinite.

We now write $\lambda(r)=\rho(A_r)$ as in Section \ref{sec:spectral}.
For real $r$ near $h$, the operator $B_{r,t}$ has a simple leading eigenvalue $e^t\lambda(r)$ with strictly positive left and right
eigenvectors.  Its remaining spectrum is contained in a strictly
smaller disk (Proposition \ref{prop:PF}).  Therefore the terms in the sum on the right of \eqref{eq:pathequal} can be written as
\[
\left\langle \xi_{r,x},B_{r,t}^{\,k}\eta_r\right\rangle =c(r,t)\left(e^t\lambda(r)\right)^k  + O_{r,t}(q_{r,t}^{\,k}),  \ \text{where}  \ 0\leq q_{r,t}<e^t\lambda(r)
\]
and $c(r,t) >0$.  It follows that the series in \eqref{eq:pathequal} converges precisely when
\begin{equation}\label{eq:pressure}
e^t\lambda(r)<1.
\end{equation}

We next consider the normalisation in \cite{CP25}. In \cite[Section 3.5]{CP25} the cost is centred before
Proposition 3.12 is applied: the centred cost is $\psi(p)=m(p)-\Lambda \ell(p)$.
Let $K(t)$ denote the exponential growth
exponent for the centred weighted path count as in \cite[Proposition 3.12]{CP25}.  The derivative
calculation in \cite[Theorem 2.4, Lemma 3.10 and the proof of
Theorem 2.6]{CP25} identifies
\[
K(0)=h, \ \ K'(0)=0 \ \ \text{ and } \ \  K''(0)=\sigma^2.
\]
Centring does not change the second derivative.  We also have
\[
\sum_{p\in\mathcal P_x} e^{-r\ell(p)+t\psi(p)} =\sum_{p\in\mathcal P_x} e^{-(r+t\Lambda)\ell(p)+tm(p)}.
\]
It follows that the abscissa of convergence of the uncentred  series is the analytic function
\[
H(t)=K(t)+t\Lambda, \ \ \  H(0)=h, \ \ \  H'(0)=\Lambda, \ \ \  H''(0)=\sigma^2.
\]
Here we have undone the length normalisation $h=1$ used for the calculations in \cite[Section 3]{CP25}.
Comparing \eqref{eq:pathequal} with  \eqref{eq:pressure} now proves  that
\[
e^t\lambda\left(H(t)\right)=1 \ \ \text{ and } \ \ H(0)=h.
\]

Finally, let $a,v$ be as in Proposition \ref{prop:PF}. Then, differentiating 
$t+\log\lambda(H(t))=0$ at $t=0$ gives
\[
1-aH'(0)=0 \ \ \text{ and } \ \  v\left(H'(0)\right)^2 -aH''(0)=0.
\]
Therefore the two sets of constants are related by
\[
\Lambda=\frac{1}{a},\ \ \ \sigma^2=\frac{v}{a^3},
\]
and $a$, $v$ as  in Proposition \ref{prop:PF}  are strictly positive.

\bibliographystyle{alpha}
\bibliography{LLT_references}

@article{athreya2012distribution,
  title={The distribution of gaps for saddle connection directions},
  author={Athreya, Jayadev S and Chaika, Jon},
  journal={Geometric and Functional Analysis},
  volume={22},
  number={6},
  pages={1491--1516},
  year={2012},
  publisher={Springer}
}

@article{CWW20,
  author  = {Chang, Kung-Ching and Wang, Xuefeng and Wu, Xie},
  title   = {On the spectral theory of positive operators and {PDE} applications},
  journal = {Discrete and Continuous Dynamical Systems},
  volume  = {40},
  number  = {6},
  pages   = {3171--3200},
  year    = {2020},
  doi     = {10.3934/dcds.2020054}
}

@article{babillot.ledrappier,
	author = {M. Babillot and F. Ledrappier},
	journal = {Ergodic theory and dynamical systems},
	pages = {17–39},
	title = {Lalley’s theorem on periodic orbits of hyperbolic flows},
	volume = {188},
	year = {1998}}

@article{schwartz.sharp,
	author = {R. Sharp and R. Schwartz},
	journal = {Communications in Mathematical Physics},
	volume ={153},
	pages ={423-430},
	title = {The correlation of length spectra of two hyperbolic surfaces},
	year = {1993}}

@article{DolgopyatHafouta2022,
  author  = {Dolgopyat, Dmitry and Hafouta, Yeor},
  title   = {Edgeworth Expansions for Independent
             Bounded Integer Valued Random Variables},
  journal = {Stochastic Processes and their Applications},
  volume  = {152},
  pages   = {486--532},
  year    = {2022},
  doi     = {10.1016/j.spa.2022.07.001}
}

@article{milnor,
author = {J. Milnor},
title = {{A note on curvature and fundamental group}},
volume = {2},
journal = {Journal of Differential Geometry},
number = {1},
publisher = {Lehigh University},
pages = {1 -- 7},
year = {1968},
doi = {10.4310/jdg/1214501132}
}

@article{Lalley1987,
  author  = {Lalley, Steven P.},
  title   = {Distribution of Periodic Orbits of Symbolic
             and Axiom A Flows},
  journal = {Advances in Applied Mathematics},
  volume  = {8},
  pages   = {154--193},
  year    = {1987}
}

@article{athreya2023counting,
  title={Counting pairs of saddle connections},
  author={Athreya, Jayadev S and Fairchild, Samantha and Masur, Howard},
  journal={Advances in Mathematics},
  volume={431},
  pages={109233},
  year={2023},
  publisher={Elsevier}
}

@article{robertson2023uniform,
  title={Uniform distribution of saddle connection lengths in all SL (2, R) orbits},
  author={Robertson, Donald and Dozier, Benjamin},
  journal={Geometriae Dedicata},
  volume={217},
  number={4},
  pages={65},
  year={2023},
  publisher={Springer}
}

@article{NevoRuhrWeiss2020,
  author  = {Nevo, Amos and R{\"u}hr, Rene
             and Weiss, Barak},
  title   = {Effective Counting on Translation Surfaces},
  journal = {Advances in Mathematics},
  volume  = {360},
  pages   = {106890},
  year    = {2020},
  doi     = {10.1016/j.aim.2019.106890}
}

@article{Lalley1989,
  author  = {Lalley, Steven P.},
  title   = {Renewal Theorems in Symbolic Dynamics,
             with Applications to Geodesic Flows,
             Noneuclidean Tessellations and
             Their Fractal Limits},
  journal = {Acta Mathematica},
  volume  = {163},
  pages   = {1--55},
  year    = {1989}
}

@article{Sharp2004,
  author  = {Sharp, Richard},
  title   = {A Local Limit Theorem for Closed
             Geodesics and Homology},
  journal = {Transactions of the American
             Mathematical Society},
  volume  = {356},
  number  = {12},
  pages   = {4897--4908},
  year    = {2004}
}

@article{CantrellPollicott2022,
  author  = {Cantrell, Stephen and Pollicott, Mark},
  title   = {Comparison Theorems for Closed Geodesics
             on Negatively Curved Surfaces},
  journal = {Groups, Geometry, and Dynamics},
  volume  = {16},
  number  = {2},
  pages   = {461--491},
  year    = {2022},
  doi     = {10.4171/GGD/671}
}

@book{EisnerFarkasHaaseNagel2015,
  author    = {Eisner, Tanja and Farkas, B{\'a}lint
               and Haase, Markus and Nagel, Rainer},
  title     = {Operator Theoretic Aspects of Ergodic Theory},
  series    = {Graduate Texts in Mathematics},
  volume    = {272},
  publisher = {Springer},
  year      = {2015},
  doi       = {10.1007/978-3-319-16898-2}
}

@misc{CP25,
  author       = {Cantrell, Stephen and Pollicott, Mark},
  title        = {Counting statistics for geodesics on flat surfaces},
  year         = {2025},
  note         = {arXiv:2503.13091.}
}

@article{ColP,
  title={The growth and distribution of large circles on translation surfaces},
  author={Colognese, Paul and Pollicott, Mark},
  journal={arXiv preprint arXiv:2107.14058},
  year={2021}
}

@article{Wright15,
  author       = {Wright, A.},
  title        = {Translation surfaces and their orbit closures: An introduction for a broad audience},
  journal      = {EMS Surv. Math. Sci.},
  volume       = {2},
  year         = {2015},
  pages        = {63-108}
}

@article{Mas90,
  author       = {Masur, H.},
  title        = {The growth rate of trajectories of a quadratic differential},
  journal      = {Ergodic Theory Dynam. Systems},
  volume       = {10},
  year         = {1990},
  pages        = {151--176}
}

@book{Simon05,
  author       = {Simon, B.},
  title        = {Trace Ideals and Their Applications},
  edition      = {2nd},
  series       = {Mathematical Surveys and Monographs},
  volume       = {120},
  publisher    = {American Mathematical Society},
  year         = {2005}
}

@misc{Peterson13,
  author       = {Peterson, J.},
  title        = {Notes on {von Neumann} Algebras},
  howpublished = {Lecture notes},
  year         = {2013}
}

@book{Kato,
  author       = {Kato, T.},
  title        = {Perturbation Theory for Linear Operators},
  publisher    = {Springer},
  year         = {1980}
}

@book{HornJohnson13,
  author       = {Horn, R. A. and Johnson, C. R.},
  title        = {Matrix Analysis},
  edition      = {2nd},
  publisher    = {Cambridge University Press},
  year         = {2013}
}

@article{Stone65,
  author       = {Stone, C.},
  title        = {A local limit theorem for nonlattice multi-dimensional distribution functions},
  journal      = {Ann. Math. Statist.},
  volume       = {36},
  year         = {1965},
  pages        = {546--551}
}

@article{weierstrass,
  author       = {Barletta, Elisabetta and Dragomir, Sorin and Esposito, Francesco},
  title        = {Cauchy formula for vector-valued holomorphic functions and the {Cauchy-Kovalevskaja} theorem},
  journal      = {Analysis and Geometry in Metric Spaces},
  volume       = {13},
  number       = {1},
  year         = {2025}
}

@article{EM01,
  author  = {Eskin, Alex and Masur, Howard},
  title   = {Asymptotic formulas on flat surfaces},
  journal = {Ergodic Theory and Dynamical Systems},
  volume  = {21},
  number  = {2},
  year    = {2001},
  pages   = {443--478},
  doi     = {10.1017/S0143385701001225}
}

@article{EMZ03,
  author  = {Eskin, Alex and Masur, Howard and Zorich, Anton},
  title   = {Moduli spaces of {Abelian} differentials: The principal boundary,
             counting problems, and the {Siegel--Veech} constants},
  journal = {Publications Math\'ematiques de l'IH\'ES},
  volume  = {97},
  year    = {2003},
  pages   = {61--179},
  doi     = {10.1007/s10240-003-0015-1}
}

\end{document}